Glasgow Mathematical Journal

**RESEARCH ARTICLE**

# Actions of nilpotent groups on nilpotent groups

Michael C. Burkhart[1,2] 

[1] University of Cambridge, Cambridge, UK
[2] Current Address: University of Chicago, Chicago, USA
Emails: mcb93@cantab.ac.uk, burkh4rt@uchicago.edu

**Abstract**
For finite nilpotent groups $J$ and $N$, suppose $J$ acts on $N$ via automorphisms. We exhibit a decomposition of the first cohomology set in terms of the first cohomologies of the Sylow $p$-subgroups of $J$ that mirrors the primary decomposition of $H^1(J, N)$ for abelian $N$. We then show that if $N \rtimes J$ acts on some non-empty set $\Omega$, where the action of $N$ is transitive and for each prime $p$ a Sylow $p$-subgroup of $J$ fixes an element of $\Omega$, then $J$ fixes an element of $\Omega$.

## 1. Introduction

Given a finite nilpotent group $J$ acting on a finite nilpotent group $N$ via automorphisms, crossed homomorphisms are maps $\varphi : J \to N$ satisfying $\varphi(jj') = \varphi(j)\varphi(j')^{j^{-1}}$ for all $j, j' \in J$. Two such maps $\varphi$ and $\varphi'$ are cohomologous if there exists $n \in N$ such that $\varphi'(j) = n^{-1}\varphi(j)n^{j^{-1}}$ for all $j \in J$; in this case, we write $\varphi \sim \varphi'$. We define the first cohomology $H^1(J, N)$ to be the pointed set $Z^1(J, N)$ of crossed homomorphisms modulo this equivalence relation where the distinguished point corresponds to the class containing the map taking each element of $J$ to the identity of $N$.

We first show that the cohomology set $H^1(J, N)$ decomposes in terms of the first cohomologies of the Sylow $p$-subgroups $J_p$ of $J$ as follows:

**Lemma 1.** *For finite nilpotent groups $J$ and $N$, suppose $J$ acts on $N$ via automorphisms. Then the map $\varphi \mapsto \times_{p \in \mathcal{D}} \varphi|_{J_p}$ for $\varphi \in H^1(J, N)$ induces an isomorphism $H^1(J, N) \cong \times_{p \in \mathcal{D}} H^1(J_p, N)^{J'_p}$ of pointed sets, where $\mathcal{D}$ denotes the shared prime divisors of $|J|$ and $|N|$, and for each $p$, $J_p$ is the Sylow $p$-subgroup of $J$ and $J'_p$ is the Hall $p'$-subgroup of $J$.*

This parallels the well-known primary decomposition of $H^1(J, N)$ for abelian $N$ (see Section 3 for details). As the bijective correspondence between $H^1(J, N)$ and the $N$-conjugacy classes of complements to $N$ in $N \rtimes J$ continues to hold for nonabelian $N$ [6, Exer. 1 in §I.5.1], Lemma 1 provides an alternate proof of a result of Losey and Stonehewer [5]:

**Proposition 2** (Losey and Stonehewer)**.** *Two nilpotent complements of a normal nilpotent subgroup in a finite group are conjugate if and only if they are locally conjugate.*

Here, two subgroups $H, H' \leq G$ are locally conjugate if a Sylow $p$-subgroup of $H$ is conjugate to a Sylow $p$-subgroup of $H'$ for each prime $p$. It also readily follows that:

**Proposition 3.** *Let G be a finite split extension over a nilpotent subgroup N such that G/N is nilpotent. If for each prime p, there is a Sylow p-subgroup S of G such that any two complements of $S \cap N$ in S are conjugate in G, then any two complements of N in G are conjugate.*

We then establish a fixed point result for nilpotent-by-nilpotent actions in the style of Glauberman:

**Theorem 4.** *For finite nilpotent groups J and N, suppose J acts on N via automorphisms and that the induced semidirect product $N \rtimes J$ acts on some non-empty set $\Omega$ where the action of N is transitive. If for each prime p, a Sylow p-subgroup of J fixes an element of $\Omega$, then J fixes an element of $\Omega$.*

Glauberman showed that this result holds whenever the orders of N and J are coprime, without any further restrictions on N or J [4, Thm. 4]. Thus, this result is only interesting when $|N|$ and $|J|$ share one or more prime divisors (i.e. when the action is *non-coprime*). Analogous results hold if N is abelian or if N is nilpotent and $N \rtimes J$ is supersoluble [2].

### 1.1. Outline

In the remainder of this section, we introduce some notation. We then prove the results in Section 2 and conclude in Section 3.

### 1.2. Notation

All groups in this note are finite. For a nilpotent group $J$, we let $J_p \in \mathrm{Syl}_p(J)$ denote its unique Sylow $p$-subgroup and $J'_p$ denote its Hall $p'$-subgroup so that $J \cong J_p \times J'_p$. We let $g^\gamma = \gamma^{-1} g \gamma$ for $g, \gamma \in G$. We otherwise use standard notation from group theory that can be found in Doerk and Hawkes [3].

For a subgroup $K \leq J$, we let $\varphi|_K$ denote the restriction of $\varphi \in Z^1(J, N)$ to $K$ and let $\mathrm{res}_K^J : H^1(J, N) \to H^1(K, N)$ be the map induced in cohomology. For $\varphi \in Z^1(K, N)$ and $j \in J$, define $\varphi^j(x) = \varphi(x^{j^{-1}})^j$. We say $\varphi$ is $J$-invariant if $\mathrm{res}_{K \cap K^j}^K \varphi \sim \mathrm{res}_{K \cap K^j}^{K^j} \varphi^j$ for all $j \in J$ and let $\mathrm{inv}_J H^1(K, N)$ be the set of $J$-invariant elements in $H^1(K, N)$. For any $\varphi \in Z^1(J, N)$, we have $\varphi^j(x) = n^{-1} \varphi(x) n^{x^{-1}}$ where $n = \varphi(j^{-1})$ so that $\varphi^j \sim \varphi$. Consequently, $\mathrm{res}_K^J H^1(J, N) \subseteq \mathrm{inv}_J H^1(K, N)$.

For nilpotent $J$ and $\varphi \in Z^1(J_p, N)$, any $j \in J$ may be written $j = j_p \times j'_p$ for $j_p \in J_p$ and $j'_p \in J'_p$ so that $\varphi^j(x) = \varphi(x^{j_p^{-1}})^{j_p j'_p} = \varphi'(x)^{j'_p}$ for some $\varphi' \sim \varphi$. It follows that $\mathrm{inv}_J H^1(J_p, N) = H^1(J_p, N)^{J'_p}$, that is, the $J$-invariant elements of $H^1(J_p, N)$ are those fixed under conjugation by $J'_p$.

To each complement $K$ of $N$ in $NJ$, we associate $\varphi_K \in Z^1(J, N)$ as follows. For $j \in J$, we have $j = n_j^{-1} k_j$ for unique $n_j \in N$ and $k_j \in K$; we then let $\varphi_K(j) = n_j = k_j j^{-1}$. Conversely, for any $\varphi \in Z^1(J, N)$, the subgroup $F(\varphi) = \{\varphi(j)j\}_{j \in J}$ complements $N$ in $NJ$. In particular, $F(\varphi_K) = K$. Furthermore, $F(\varphi)$ and $F(\varphi')$ are $N$-conjugate in $NJ$ if and only if $\varphi \sim \varphi'$ so that $F$ induces a correspondence between $H^1(J, N)$ and the $N$-conjugacy classes of complements to $N$ in $NJ$. See Serre [6, Ch. I §5] for further details on nonabelian group cohomology.

## 2. Proofs of results

We begin by establishing Lemma 1.

*Proof of Lemma 1.* As $N$ is nilpotent, the natural projection maps $N \to N_p$ induce an isomorphism:
$$H^1(J, N) \cong \times_{p \in \mathcal{D}} H^1(J, N_p), \tag{1}$$
where terms $p \notin \mathcal{D}$ drop by the Schur–Zassenhaus theorem [3, Thm. A.11.3]. Thus, we may focus our attention on $H^1(J, N_q)$ for some prime $q$. If $J$ is also a $q$-group, we are done. Otherwise, we may consider the inflation-restriction exact sequence [6, Sec. I.5.8]:

$$1 \to H^1(J'_q, N_q^{J_q}) \to H^1(J, N_q) \xrightarrow{\text{res}^J_{J_q}} H^1(J_q, N_q)^{J'_q}.$$

As $H^1(J'_q, N_q^{J_q})$ is trivial, $\text{res}^J_{J_q}$ is injective. For any $\varphi \in H^1(J_q, N_q)^{J'_q}$, we may define $\tilde{\varphi}$ in terms of a representative crossed homomorphism as $\tilde{\varphi}(j'j) = \varphi(j)$ for $j \in J_q$ and $j' \in J'_q$. It is straightforward to verify that $\tilde{\varphi} \in Z^1(J, N_q)$ and $\tilde{\varphi}|_{J_q} \sim \varphi$. Thus, $\text{res}^J_{J_q}$ is also surjective and thus an isomorphism. Let $\nu : H^1(J_q, N_q) \to H^1(J_q, N)$ denote the map induced by inclusion. From the decomposition (1), $H^1(J_q, N_q) \cong H^1(J_q, N)$, and as

$$H^1(J_q, N_q) \xrightarrow{\nu} H^1(J_q, N) \to H^1(J_q, N'_q)$$

is exact [6, Prop. I.38] where $H^1(J_q, N'_q)$ is trivial, it follows that $\nu$ is surjective and hence an isomorphism. As $\varphi \in H^1(J_p, N_p)^{J'_p}$ if and only if $\nu(\varphi) \in H^1(J_p, N)^{J'_p}$, it follows that $\Phi : H^1(J, N) \to \times_{p \in \mathcal{D}} H^1(J_p, N)^{J'_p}$ given by the composition $\varphi \mapsto \times_{p \in \mathcal{D}} \varphi|_{J_p}$ induces the desired isomorphism:

$$H^1(J, N) \cong \times_{p \in \mathcal{D}} H^1(J_p, N_p)^{J'_p} \cong \times_{p \in \mathcal{D}} H^1(J_p, N)^{J'_p}. \tag{2}$$

□

We now show how Propositions 2 and 3 follow from Lemma 1.

*Proof of Proposition* 2. Suppose $J$ and $J'$ each complement $N \triangleleft G$ as described in the hypotheses of the proposition. Let $\varphi$ be a crossed homomorphism representing $J'$ in $H^1(J, N)$. By hypothesis, $\varphi|_{J_p} \sim 1|_{J_p}$ for every prime $p$, where $1 \in Z^1(J, N)$ represents the distinguished point. Lemma 1 implies $\varphi \sim 1$ so that $J$ and $J'$ are conjugate. □

*Proof of Proposition* 3. Suppose $G \cong N \rtimes J$ satisfies the hypotheses of the proposition. Fix a prime $p$. Without loss, we may suppose any two complements of $N_p$ in $S = J_p N_p$ are conjugate in $G$. If $J'_p$ is such a complement, then $J'_p = (J_p)^g$ for some $g \in G$ so that $J'_p = (J_p)^{jn} = (J_p)^n$ for some $j \in J$ and $n \in N$. In particular, $J'_p$ is conjugate to $J_p$ in $J_p N$. Thus, $H^1(J_p, N)$ is trivial. As the choice of prime $p$ was arbitrary, Lemma 1 implies that $H^1(J, N)$ is also trivial, allowing us to conclude. □

To prove Theorem 4, we also require:

**Proposition 5.** *Let $H$ be a subgroup of $G \cong N \rtimes J$ where $N$ and $J$ are nilpotent. If for each prime $p$, $H$ contains a conjugate of some $J_p \in \text{Syl}_p(J)$, then $H$ contains a conjugate of $J$.*

*Proof of Proposition* 5. It follows from the hypotheses that $H$ supplements $N$ in $G$. We induct on the order of $G$. If $H$ is a $p$-group or all of $G$, the result is immediate. If multiple primes divide $|N|$, we have the nontrivial decomposition $N \cong N_p \times N'_p$ for some prime $p$. Induction in $G/N_p$ implies $J^{n_0} \leq HN_p$ for some $n_0 \in N'_p$. Induction in $G/N'_p$ implies $J^{n_1} \leq HN'_p$ for some $n_1 \in N_p$. Thus, $J^{n_0 n_1} \leq HN_p \cap HN'_p = H$, where the last equality proceeds from the following argument of Losey and Stonehewer [5]. Suppose $g \in HN_p \cap HN'_p$ so that $g = h_0 n_0 = h_1 n_1$ for some $h_0, h_1 \in H$, $n_0 \in N_p$ and $n_1 \in N'_p$. Then $(h_1)^{-1} h_0 = n_1 (n_0)^{-1} \in H$. As $n_0$ and $n_1$ commute and have coprime orders, it follows that $n_0, n_1 \in H$ so $g \in H$.

We now proceed under the assumption that $N = N_q$ for some prime $q$. Upon switching to a conjugate of $H$ if necessary, we may suppose that $J_q \leq H$. Let $Z$ denote the center of $N$. If $Z \cap H$ were nontrivial, then induction in $G/(Z \cap H)$ would allow us to conclude. Otherwise, in $G/Z$, induction implies that $J^g Z/Z \leq ZH/Z$ for some $g \in G$. Let $\psi : H \to HZ/Z$ denote the isomorphism between $H$ and $HZ/Z$. Then $K = \psi^{-1}(J^g Z/Z)$ complements $N \cap H$ in $H$ and $N$ in $G$. Let $\varphi \in Z^1(K, N)$ correspond to $J$. Then $\varphi|_{K_q}$ corresponds to $J_q$ where $[\varphi|_{K_q}] \in H^1(K_q, H \cap N)^{K'_q} \cong H^1(K, H \cap N)$. In particular, there exists a complement, say $L$, to $H \cap N$ in $H$ that contains $J_q$. $L$ will also complement $N$ in $G$, and as $\text{Syl}_p(L) \subseteq \text{Syl}_p(G)$ for all primes $p \neq q$, we may apply Proposition 2 to conclude that $J^{g'} = L \leq H$ for some $g' \in G$. □

With this, we are prepared to prove Theorem 4.

*Proof of Theorem* 4. Given $J$, $N$, and $\Omega$ as described in the hypotheses of the theorem, let $G = N \rtimes J$ denote the induced semidirect product and consider the stabilizer subgroup $G_\alpha$ for some $\alpha \in \Omega$. As $N$ acts transitively, $G = NG_\alpha$. For each prime $p$, the hypotheses of the theorem imply $(J_p)^{n_p} \leq G_\alpha$ for some $J_p \in \mathrm{Syl}_p(J)$ and $n_p \in N$ so that Proposition 5 allows us to conclude $J^g \leq G_\alpha$ for some $g \in G$. It follows that $J$ fixes $g \cdot \alpha$. □

## 3. Conclusion

We conclude with a brief discussion of analogous results in the abelian case. For arbitrary $J$ acting on abelian $N$, the restriction map $\mathrm{res}^J_{J_p} : H^1(J, N)_{(p)} \xrightarrow{\cong} \mathrm{inv}_J H^1(J_p, N)$ induces an isomorphism for each prime $p$, where $H^1(J, N)_{(p)}$ is the $p$-primary component of $H^1(J, N)$ and $J_p \in \mathrm{Syl}_p(J)$. Consequently, it follows from the primary decomposition of $H^1(J, N)$ that [1, Thm. III.10.3]:

$$H^1(J, N) \cong \oplus_{p \in \mathcal{D}} \mathrm{inv}_J H^1(J_p, N). \tag{3}$$

Furthermore, for abelian $N$, suppose $G = N \rtimes J$ acts on some non-empty set $\Omega$, where the action of $N$ is transitive, and for each prime $p$ a Sylow $p$-subgroup of $J$ fixes an element of $\Omega$. Then, for arbitrary $\alpha \in \Omega$, the stabilizer $G_\alpha$ splits over $G_\alpha \cap N$ by Gaschütz's theorem [3, Thm. A.11.2] and is locally conjugate and thus conjugate to $J$ by an argument analogous to the proof of Proposition 5. In particular, $J$ fixes an element of $\Omega$. In this note, we find that the decomposition (3) and fixed point result continue to hold for nilpotent $N$ if $J$ is also nilpotent.

**Acknowledgments.** The author thanks the editor C. M. Roney-Dougal and an anonymous reviewer for detailed feedback which considerably improved the manuscript.

**Competing interests.** The author declares none.

## References

[1] K. S. Brown, Cohomology of groups, (Springer, New York, 1982).
[2] M. C. Burkhart, Fixed point conditions for non-coprime actions, *Proc. Roy. Soc. Edinb. Sect. A* (in press).
[3] K. Doerk and T. Hawkes, Finite soluble groups, (de Gruyter, Berlin, 1992).
[4] G. Glauberman, Fixed points in groups with operator groups, *Math. Zeitschr.* **84** (1964), 120–125.
[5] G. O. Losey and S. E. Stonehewer, Local conjugacy in finite soluble groups, *Quart. J. Math. Oxford (2)* **30** (1979), 183–190.
[6] J.-P. Serre, *Galois Cohomology* (Springer, Berlin, 2002).